\documentclass{amsart}
\usepackage{a4,amsmath,amsfonts,bbold}

\newtheorem{theorem}{Theorem}
\newtheorem{definition}{Definition}
\newtheorem{lemma}{Lemma}

\DeclareMathOperator{\littleOOp}{o}
\DeclareMathOperator{\outOp}{out}
\DeclareMathOperator{\optOp}{opt}
\DeclareMathOperator{\lengthOp}{length}

\DeclareMathOperator{\reachOp}{R}
\DeclareMathOperator{\distOp}{dist}
\DeclareMathOperator{\goodOp}{G}

\newcommand{\R}{\mathbb{R}}
\newcommand{\N}{\mathbb{N}}
\newcommand{\setdef}[2]{\{{#1}\,:\,{#2}\}}
\newcommand{\setdefBig}[2]{\big\{{#1}\,:\,{#2}\big\}}
\newcommand{\out}[1]{\outOp({#1})}
\newcommand{\expect}[1]{\mathbb{E}[{#1}]}
\newcommand{\prob}[1]{\mathbb{P}[{#1}]}
\newcommand{\probBig}[1]{\mathbb{P}\big[{#1}\big]}

\newcommand{\probCondBig}[2]{\mathbb{P}\big[{#1}\,|\,{#2}\big]}
\newcommand{\vopt}{v_{\optOp}}
\newcommand{\length}[1]{\lengthOp({#1})}
\newcommand{\reach}[2]{\reachOp_{#1}({#2})}
\newcommand{\bdreach}[2]{\partial\reachOp_{#1}({#2})}
\newcommand{\good}[2]{\goodOp({#1},{#2})}
\newcommand{\ddist}[2]{\distOp^{\rightarrow}({#1},{#2})}
\newcommand{\dist}[2]{\distOp({#1},{#2})}
\newcommand{\littleO}[1]{\littleOOp({#1})}
\newcommand{\lofac}[2]{{#1}^{\underline{#2}}}
\newcommand{\allzero}{\mathbb{0}}
\newcommand{\evec}[1]{\mathbb{e}_{#1}}

\newcommand{\RF}{\mbox{\sc Random-Facet}}
\newcommand{\RE}{\mbox{\sc Random-Edge}}
\newcommand{\GE}{\mbox{\sc Greatest-Decrease}}

\title{Two New Bounds for the Random-Edge Simplex Algorithm}
\author{Bernd G\"artner}

\address{
  Institut f\"ur Theoretische Informatik\\
  ETH Z\"urich\\
  Haldeneggsteig 4\\
  CH-8092 Z\"urich
}
\email{gaertner@inf.ethz.ch}
\thanks{The first author acknowledges support from the Swiss Science 
Foundation (SNF), Project No.\ 200021-100316/1. The second author has been supported by the DFG Research Center \textsc{Matheon} and by the DFG Research Group \emph{Algorithms, Structure, Randomness}.}
\author{Volker Kaibel}
\address{
  Zuse-Institut Berlin\\
  Takustr. 7\\
  D-14195 Berlin
}
\email{kaibel@zib.de}

\begin{document}

\begin{abstract}
We prove that the \RE\ simplex algorithm requires an
expected number of at most $13n/\sqrt{d}$ pivot steps
on any simple $d$-polytope with $n$ vertices. This is the first
nontrivial upper bound for general polytopes. We also describe a
refined analysis that potentially yields much better bounds for 
specific classes of polytopes. As one application, we show that
for combinatorial $d$-cubes, the trivial upper bound of $2^d$ on
the performance of \RE\ can asymptotically be improved by any 
desired polynomial factor in $d$.
\end{abstract}

\maketitle

\section{Introduction}

Dantzig's \emph{simplex method} \cite{d-lpe-63} is a widely used tool 
for solving linear programs (LP). The feasible region 
of an LP is a polyhedron; any algorithm implementing the simplex method
traverses a sequence of vertices, such that (i) consecutive vertices 
are equal (the \emph{degenerate} case) or connected by a polyhedron edge, 
and (ii) the objective function strictly improves along any traversed edge. 
In both theory and practice, we may
assume that some initial vertex is available, and that the
optimal solution to the LP is attained at a vertex, if there
is an optimum at all. It follows that if the algorithm does not 
cycle, it will eventually find an optimal solution, or discover
that the problem is unbounded (see e.g.\  Chv{\'a}tal's book \cite{c-lp-83}
for a comprehensive introduction to the simplex method). 

For 
most (complexity-)theoretic investigations, one can safely assume
that the LP's that are considered are bounded as well as both primally
and dually non-degenerate \cite{s-tlip-86}. 
Thus, we will only deal with \emph{simple
polytopes}, i.e., bounded $d$-di\-men\-sional polyhedra, where at each
vertex exactly $d$ facets meet, and with objective functions that are
non-constant along any edge of the polytope.

The distinguishing feature of each simplex-algorithm is the \emph{pivot
rule} according to which the next vertex in the sequence is
selected in case there is a choice. Many popular pivot rules are
efficient in practice, meaning that they induce a short vertex
sequence in typical applications. The situation in theory is
in sharp contrast to this: Among  most of the \emph{deterministic}
pivot rules proposed in the literature (including the ones widely
used in practice), the simplex algorithm is forced to traverse 
exponentially (in the number of variables and constraints of the LP) 
many vertices in the worst case. It is open whether there is a pivot
rule that always induces a sequence of polynomial length.

To explain simplex's excellent behavior in practice, the tools of
\emph{average case analysis} \cite{b-smpa-87} and \emph{smoothed analysis} 
\cite{st-saa-91} have been devised, and to conquer the worst case bounds,
research has turned to \emph{randomized} pivot rules. Indeed, Kalai 
\cite{k-srsa-92,k-lpsasp-97} as well as Matou\v{s}ek, Sharir and Welzl 
\cite{msw-sblp-96} could
prove that the expected number of steps taken by the \RF\ pivot rule is
only \emph{subexponential} in the worst case. These results hold under
our above assumption that the feasible region of the LP is a simple and 
full-dimensional polytope.

Much less is known about another (actually, the most natural)
randomized pivot rule: choose the next vertex in the sequence
uniformly at random among the neighbors of the current vertex with
better objective function value. This rule is called \RE, and unlike
\RF, it has no recursive structure to peg an analysis to. Nontrivial
upper bounds on its expected number of pivot steps on general
polytopes do not exist. Results are known for 3-polytopes
\cite{bdfru-wcrtr-95,kmsz-sadt-04}, $d$-polytopes with $d+2$
facets~\cite{gstvw-olp-03}, and for linear assignment
problems~\cite{Tov86}.
 Only recently, Pemantle and Balogh solved the
long standing problem of finding a tight bound for the expected
performance of \RE\ on the $d$-dimensional \emph{Klee-Minty} cube
\cite{bp-recmwls-04}. This polytope is the 'mother' of many worst-case
inputs for deterministic pivot rules \cite{km-hgisa-72,az-dpmsp-99}.

None of the existing results exclude the possibility of both \RF\ and
\RE\ being the desired (expected) polynomial-time pivot rules.
In the more general and well-studied setting of 
\emph{abstract objective functions} on polytopes 
\cite{as-lmpap-76,w-usf-85,w-cunsp-88,k-lpsasp-97}, 
superpolynomial lower bounds are known 
for both rules, where the construction for \RE\ \cite{ms-reeac-04} 
is very recent and much more involved than the one for 
\RF\ \cite{m-lbsoa-94}. Both approaches
inherently use objective functions (on cubes) that are not linearly 
induced.

In this paper, we derive the first nontrivial upper bound for the
expected performance of \RE\ on simple polytopes, with edge
orientations induced by abstract objective functions. 
Even when we
restrict to linear objective functions on combinatorial cubes, the
result is new. The general bound itself is rather weak and also
achieved for example by the deterministic 
\GE\ rule. The emphasis here is on the fact
that we are able to make progress at all, given that \RE\ has turned
out to be very difficult to attack in the past. Also, our new bound
separates \RE\ from many deterministic rules (for example,
\emph{Dantzig's rule}, \emph{Bland's rule}, or the \emph{shadow vertex
  rule}) that may visit all vertices in the worst case
\cite{az-dpmsp-99}.

In a second part, we refine the analysis, with the goal of obtaining
better bounds for specific classes of polytopes. Roughly speaking,
these are polytopes with large and regular local neighborhoods. 
Our prime example
is the class of combinatorial cubes, for which we improve the
general upper bound by any desired polynomial factor in the dimension.
As before, this also works for abstract objective functions and thus
complements the recent lower bound of Matou\v{s}ek and Szab\'{o}
\cite{ms-reeac-04} with a first nontrivial upper bound.

\section{A Bound for General Polytopes}
\label{sec:general}

Throughout this section, $P$ is a $d$-dimensional simple polytope
with a set $V$ of $n$ vertices.
A directed graph $D=(V,A)$ is called an \emph{acyclic unique sink 
orientation} (AUSO) of~$P$ if
\begin{itemize}
\item[(i)] its underlying 
undirected graph is the vertex-edge graph of $P$, 
\item[(ii)] $D$ is acyclic, and
\item[(iii)] any subgraph of $D$ induced by the vertices of a
nonempty face of $P$ has a unique sink. 
\end{itemize}
Any linear function $\varphi:V\rightarrow\R$ that is 
\emph{generic} (non-constant on 
edges of~$P$) induces an AUSO in a natural way: there is a directed
edge $v\rightarrow w$ between adjacent vertices if and only
if $\varphi(v)>\varphi(w)$. The global sink of the AUSO is the unique 
vertex that minimizes $\varphi$ over $P$. If $\varphi$ is any generic (not 
necessarily linear) function inducing an AUSO that way, $\varphi$
is called an \emph{abstract objective function}. For a
given AUSO $D$ of~$P$, any function $\varphi$ that maps vertices 
to their ranks w.r.t.\ a fixed topological sorting of $D$ is
an abstract objective function that induces $D$. In general,
$D$ need not be induced by a linear function, for example if
$D$ fails to satisfy the necessary \emph{Holt-Klee} condition
for linear realizability \cite{hk-psmfsc-98}.
For the remainder of this section, we fix an AUSO $D$ of $P$, 
an abstract objective function $\varphi$ that induces $D$ and 
some vertex $s\in V$.

Let~$\pi$ be the random variable defined as 
the directed path  in~$D$,
starting at~$s$ and ending at the sink~$\vopt$ of~$D$, induced by the
\RE\ pivot rule. From each visited vertex $v\neq\vopt$, $\pi$ proceeds
to a neighbor $w$ of $v$, along an outgoing edge chosen uniformly at
random from all outgoing edges.

For each $v\in V$, denote by 
$$
\out{v}\ :=\ \setdef{w\in V}{(v,w)\in A}
$$
the set of all smaller (w.r.t.~$\varphi$) neighbors of~$v$. If
$|\out{v}|=k$, then~$v$ is called a \emph{$k$-vertex}. We denote by
$V_k$ the set of all $k$-vertices. 

For every vertex~$v\ne\vopt$ on the path~$\pi$ let~$v'$ be its
successor on~$\pi$. We denote by
$$
S(v)\ :=\ \setdef{w\in\out{v}}{\varphi(v')<\varphi(w)}\ ,
$$
the set of neighbors of~$v$ that are 'skipped' by~$\pi$ at the step
from~$v$ to~$v'$.
For every $0\le k\le d$ let 
$$
\eta_k(\pi)\ :=\ \big|\setdef{v\in\pi\cap V_k}{|S(v)|\ge\lfloor\tfrac{|\out{v}|}{2}\rfloor}\big|\
$$
be the number of $k$-vertices on~$\pi$, where~$\pi$ skips at least
$\lfloor\tfrac{k}{2}\rfloor$ neighbors. (Here, as in the following, we
write, depending on the context, '$\pi$' for the set of
vertices on the path~$\pi$.)

If we  denote by $n_k(\pi)$ the total number of $k$-vertices on the
path~$\pi$, then we obtain
\begin{equation}
  \label{eq:eta-n}
  \expect{\eta_k(\pi)}\ \ge\ \tfrac{1}{2}\expect{n_k(\pi)}\ .
\end{equation}
Indeed, we have
$$
\expect{\eta_k(\pi)}\ =\ 
\sum_{v\in V_k}\prob{v\in\pi\text{ and }|S(v)|\ge\lfloor\tfrac{k}{2}\rfloor}
$$
and
$$
\expect{n_k(\pi)}\ =\ 
\sum_{v\in V_k}\prob{v\in\pi}\ .
$$
The claim then follows from
$$
\probCondBig{|S(v)|\ge\lfloor\tfrac{|\out{v}|}{2}\rfloor}{v\in\pi}\ \ge\ \tfrac{1}{2}\ .
$$

Due to $\varphi(v)>\varphi(w)>\varphi(v')$ for all $w\in S(v)$, the
sets~$S(v)$ are pairwise disjoint. Thus, we obtain (exploiting the
linearity of expectation) for the number $\length{\pi}$ of vertices on $\pi$
$$
\expect{\length{\pi}}
\ \le\ 
n-\sum_{k=0}^d \expect{\eta_k(\pi)}\lfloor\tfrac{k}{2}\rfloor
\ \le\ 
n-\sum_{k=0}^d \tfrac{1}{2}\lfloor\tfrac{k}{2}\rfloor\expect{n_k(\pi)}
$$
(where we used~(\ref{eq:eta-n}) for the second inequality).
Clearly, we have $\expect{\length{\pi}}=\sum_{k=0}^d
\expect{n_k(\pi)}$.  Therefore, we obtain (note
$\tfrac{1}{2}\lfloor\tfrac{k}{2}\rfloor\ge\tfrac{k-1}{4}$)
\begin{equation}
  \label{eq:pimin}
  \expect{\length{\pi}}\ \le\ 
  \min
  \big\{
  \sum_{k=0}^d \expect{n_k(\pi)}, 
               n-\sum_{k=0}^d \tfrac{k-1}{4}\expect{n_k(\pi)}
  \big\}    
  \ .
\end{equation}

If~$h_k$ denotes the total number of $k$-vertices in~$V$, then we clearly have 
$0\le \expect{n_k(\pi)}\le h_k$. Thus, (\ref{eq:pimin}) yields
\begin{equation}
  \label{eq:maxmin}
  \expect{\length{\pi}}\ \le\ 
  \max
    \setdefBig%
      {\min\{\sum_{k=0}^d x_k, n-\sum_{k=0}^d\tfrac{k-1}{4}x_k\}}%
      {0\le x_k\le h_k\text{ for all~$k$}}
  \ .
\end{equation}

In~(\ref{eq:maxmin}), the maximum must be attained by some
$x\in\R^{d+1}$ for which the minimum is attained by both $\sum x_k$
and $n-\sum\tfrac{k-1}{4}x_k$. Indeed, if $\sum
x_k<n-\sum\tfrac{k-1}{4}x_k$ then not all~$x_k$ can be at their
respective upper bounds~$h_k$ (since $n=\sum h_k$), thus one of them
can slightly be increased in order to increase the minimum. If $\sum
x_k>n-\sum\tfrac{k-1}{4}x_k$ then not all~$x_k$ can be zero (since 
this would yield $0>n$), so one of
them can be decreased in order to increase the minimum. Thus we
conclude
\begin{equation}
  \label{eq:max}
  \expect{\length{\pi}}\ \le\ 
  \max
    \setdefBig%
      {\sum_{k=0}^d x_k}
      {\sum_{k=0}^d\tfrac{k+3}{4}x_k=n,\, 0\le x_k\le h_k\text{ for all~$k$}}
  \ .
\end{equation}

By (weak) linear programming duality (and exploiting
$n=\sum_{k=0}^dh_k$ once more), we can derive 
from~(\ref{eq:max}) the estimate
\begin{equation}
  \label{eq:dual}
  \expect{\length{\pi}}
  \ \le\ 
  \sum_{k=0}^d h_k\cdot\max\{y,1-\tfrac{k-1}{4}y\}
\end{equation}
for every $y\in\R$.

In the sequel, we need two important results from the theory of
convex polytopes. The parameters~$h_k$ are independent of
the actual acyclic unique sink orientation of the polytope. The
\emph{$h$-vector} formed by them is a linear transformation of the
\emph{$f$-vector} of the polytope, storing for each~$i$ the number of
$i$-dimensional faces of the polytope.

The first classical result we need are the \emph{Dehn-Sommerville equations}
\begin{equation}
  \label{eq:DS}
  h_k\ =\ h_{d-k}\qquad\text{for all } 0\le k\le d
\end{equation}
(see \cite[Sect.~8.3]{Z-lp-94}). The second one is the 
\emph{unimodality of the $h$-vector}:
\begin{equation}
  \label{eq:unimod}
  h_0\ \le\ h_1\ \le \ \dots\ \le\ h_{\lfloor d/2\rfloor}
\end{equation}
The latter is equivalent to the \emph{nonnegativity of the
  $g$-vector}, which is one of the hard parts of the \emph{$g$-theorem
  for simplicial polytopes}, see \cite[Sect.~8.6]{Z-lp-94}.

From~(\ref{eq:DS}) and~(\ref{eq:unimod}) we can derive 
$$
n\ =\ \sum_{k=0}^d h_k\ \ge\ \big(d-8\sqrt{d}\big)h_{\lfloor4\sqrt{d}\rfloor}\ ,
$$
which yields (for $d > 64$)
\begin{equation}
  \label{eq:hsqrtd}
  h_{\lfloor4\sqrt{d}\rfloor}\ \le\ \frac{n}{d-8\sqrt{d}}\ .
\end{equation}

Now we choose $y:=1/{\sqrt{d}}$ in~(\ref{eq:dual}). We have
$$
\frac{1}{\sqrt{d}}\ge 1-\frac{k-1}{4\sqrt{d}}
\ \Leftrightarrow\
k\ge 4\sqrt{d}-3\ .
$$
Thus, (\ref{eq:dual}) (with $y=1/{\sqrt{d}}$) gives
\begin{equation}
  \label{eq:general:final}
  \expect{\length{\pi}}\ \le\
  \sum_{k=0}^{\lfloor 4\sqrt{d}-3\rfloor}h_k\big(1-\frac{k-1}{4\sqrt{d}}\big)
  \ +\
  \sum_{k=\lfloor 4\sqrt{d}-3\rfloor +1}^d \frac{h_k}{\sqrt{d}}\ .
\end{equation}

By the unimodality of the $h$-vector and~(\ref{eq:hsqrtd}), the first
sum in~(\ref{eq:general:final}) can be estimated by
$$
4\sqrt{d}\cdot h_{\lfloor 4\sqrt{d}\rfloor}\ \le\ \frac{4n}{\sqrt{d}-8}
\le \frac{12n}{\sqrt{d}}, \quad d\geq 144\ .$$ 

Clearly, the second sum in~(\ref{eq:general:final}) is bounded by
$n/\sqrt{d}$. The resulting total bound of $13n/\sqrt{d}$ also holds
for $d<144$, because $n$ is a trivial upper bound. Thus we have proved 
the following result. 

\begin{theorem}
 The expected number of vertices visited by the \RE\ 
  simplex-algorithm on a $d$-dimensional $simple$
 polytope with~$n\geq
  d+1$ vertices, equipped with an abstract (in particular: a linear)
  objective function is bounded by
  $$
  13\cdot\frac{n}{\sqrt{d}}\ .
  $$
\end{theorem}

A similar analysis reveals that the running-time 
for the \GE-rule is bounded by 
  $$
  C'\cdot\frac{n}{\sqrt{d}}\ .
  $$
In each step, this rule selects the neighboring vertex with
\emph{smallest} $\varphi$-value, thus skipping all other neighbors
of the current vertex $v$.

For general simple polytopes, our analysis of the bound for \RE\ 
stated in~(\ref{eq:maxmin}) is essentially best possible. This can be
 seen through the examples of duals of stacked simplicial polytopes
(see, e.g., \cite{BL93}), which are simple $d$-polytopes with
$n$~vertices, $h_0=h_d=1$, and $h_k=\tfrac{n-2}{d-1}$ for all $1\le
k\le d-1$.

\section{A Bound for Cubes}
\label{sec:cubes}

The core argument of the analysis presented in
Section~\ref{sec:general} is the following: For every vertex on the
\RE\ path~$\pi$ with out-degree~$k$ we know that~$\pi$ skips
(in expectation) $k/2$ vertices \emph{in the single step from~$v$ to
  its successor}. 
We then exploited the Dehn-Sommervile equations as
well as the unimodality of the $h$-vector in order to argue that many
vertices on~$\pi$ must have large out-degree -- unless~$\pi$ is
'short' anyway.

For the $d$-dimensional cube, we have much more information on the
$h$-vector: $h_k=\binom{d}{k}$ for every~$k$. Thus, 'most' vertices
have out-degree roughly $d/2$ in case of cubes. We will exploit this
stronger knowledge in a sharper analysis for cubes, which relies on
studying larger structures around vertices than just their
out-neighbors. We actually do the analysis for general polytopes 
and obtain a bound on the expected path length in terms of two
specific quantities. Later we bound these quantities for the case 
of cubes.

\subsection{The General Approach}
\label{subsec:cubes:gen}

Within this subsection, (as in Section~\ref{sec:general}), let~$P$ be
a $d$-dimensional simple polytope with~$n$ vertices $V$, $D=(V,A)$ an
AUSO of~$P$, $\varphi:V\rightarrow\R$
an abstract objective function inducing~$D$, and $s\in V$ a fixed
vertex. We denote by
$\ddist{v}{w}$ the length (number of arcs) of a shortest directed path
from~$v$ to~$w$ ($\ddist{v}{w}$ may be $\infty$ if there is no such
path).

\begin{definition}[$t$-reach]
  Let $t,k\in\N$ and $v\in V$.
  \begin{enumerate}
  \item We call 
    $$
    \reach{t}{v}\ :=\ \setdef{w\in V}{\ddist{v}{w}\le t}
    $$
    the \emph{$t$-reach} of~$v$.  The \emph{boundary}
    of~$\reach{t}{v}$, denoted by~$\bdreach{t}{v}$, is the set of
    all $w\in\reach{t}{v}$, for which there is a directed
    (not necessarily shortest) path of length precisely~$t$ from~$v$
    to~$w$.
  \item The $t$-reach~$\reach{t}{v}$ is \emph{$k$-good} if
    $$
    |\out{w}|\ \ge\ k
    $$
    holds for all $w\in\reach{t}{v}$ with $\ddist{v}{w}\le t-1$.
  \item A vertex~$v$ is \emph{$(t,k)$-good} if its $t$-reach is
    $k$-good. The set of all $(t,k)$-good vertices is denoted by
    $\good{t}{k}$.
  \end{enumerate}
\end{definition}

In particular, if $v$ is $(t,k)$-good, the optimal vertex $\vopt$ 
may occur in the boundary of $\reach{t}{v}$, but not in its interior.
For $t,k\in\N$, we define
$$
g(t,k)\ :=\ \min\setdef{|\bdreach{t}{v}|}{v\in\good{t}{k}}\ .
$$

For every vertex~$v\in V$, and some $t\in\N$, denote by (the random
variable) $w_t(v)$ the vertex that is reached by the \RE\
simplex-algorithm, started at~$v$, after~$t$ steps (let
$w_t(v):=\vopt$ in case the sink is reached before step~$t$).
Generalizing the notion from Section \ref{sec:general}, we denote by
$$
\Tilde{S}_t(v)\ :=\ \setdef{u\in\reach{t}{v}}{\varphi(u)>\varphi(w_t(v))}
$$
the set of vertices in $\reach{t}{v}$
left behind while walking from~$v$ to~$\reach{t}{v}$.

\begin{lemma}
\label{lem:skipreach}
  For every $t,k\in\N$ and $v\in\good{t}{k}$, we have
  $$
  \probBig{|\Tilde{S}_t(v)|\ge\tfrac{g(t,k)}{2}}\ \ge\
  \frac{g(t,k)}{2d^t}\ .
  $$
\end{lemma}

\begin{proof}
  Let $\bdreach{t}{v}=\{u_1,\dots,u_q\}$ with $\varphi(u_1)>\dots
\varphi(u_q)$. By construction, there is some~$i^{\star}$ with
  $w_t(v)=u_{i^{\star}}$. Since the outdegree at every vertex is at
  most~$d$, we have
  $$
  \prob{i^{\star}=i}\ \ge\ \frac{1}{d^t}
  $$
  for every $1\le i\le q$. Therefore, 
  $$
  \prob{i^{\star}> q/2}\ \ge\ \frac{q}{2d^t}
  $$
  holds. Since $q\ge g(t,k)$ holds and because
  $i^{\star}>g(t,k)/2$ implies
  $|\Tilde{S}_t(v)|\ge g(t,k)/{2}$, the claim 
  follows.

\end{proof}

Now let us consider the path~$\pi$ followed by the \RE\ 
simplex-algorithm started at~$s$ (ending in $v_{\optOp}$).  For
$t,k\in\N$ with $t\ge 2$ and $k\ge 1$, we subdivide~$\pi$ into
subpaths with the property that every subpath 
either has length one and starts at a non-$(t,k)$-good
vertex or it has length~$t$ (a \emph{long subpath}) and starts at a
$(t,k)$-good vertex. (Such a partitioning is clearly possible.)

Let $n_{t,k}(\pi)$ be the number of long subpaths in our partitioning. We
denote the pairs of start and end vertices of these long paths by
$(x_1,y_1)$, \dots $(x_{n_{t,k}(\pi)},y_{n_{t,k}(\pi)})$.  Let
$$
S_t(x_i)\ :=\ \setdef{u\in\reach{t}{x_i}}{\varphi(u)>\varphi(y_i)}\
$$
and define
$$
\eta_{t,k}(\pi)\ :=\ 
\big|\setdefBig{i\in\{1,\dots,n_{t,k}(\pi)\}}{|S_t(x_i)|\ge \tfrac{g(t,k)}{2}}\big|
$$
to be the number of those long subpaths which leave behind
at least
$\tfrac{g(t,k)}{2}$ vertices from $\reach{t}{x_i}$. 

Using Lemma~\ref{lem:skipreach} (note that $S_t(x_i)$, conditioned on
the event that~$x_i$ is the start vertex of a long subpath in the
partitioning of~$\pi$, has the same distribution as
$\Tilde{S}_t(x_i)$), we can deduce, similarly to our derivation
of~(\ref{eq:eta-n}), the following:
\begin{equation}
  \label{eq:eta-n-reach}
  \expect{\eta_{t,k}(\pi)} \ \ge\ \frac{g(t,k)}{2d^t}\expect{n_{t,k}(\pi)}.
\end{equation}

Also here, the sets $S_t(x_i)$ (for $1\le i\le n_{t,k}(\pi)$) are
pairwise disjoint. Thus, for each long subpath (consisting of~$t$
arcs) starting at some~$x_i$ with $|S_t(x_i)|\ge g(t,k)/2$ we can
count at least $g(t,k)/2-t$ vertices that are not visited by~$\pi$.
Therefore, we can conclude
\begin{equation*}
  \expect{\length{\pi}}\ \le\ n-\big(\tfrac{g(t,k)}{2}-t\big)\expect{\eta_{t,k}(\pi)}\ .
\end{equation*}
Using~(\ref{eq:eta-n-reach}) and defining
$$
\Tilde{g}(t,k)\ :=\ \big(\tfrac{g(t,k)}{2}-t\big)\tfrac{g(t,k)}{2d^t}\ ,
$$
this yields
\begin{equation}
  \label{eq:reachbound1}
  \expect{\length{\pi}}\ \le\ n-\Tilde{g}(t,k)\expect{n_{t,k}(\pi)}\ .
\end{equation}

On the other hand, denote by
\begin{equation}
  \label{eq:def-f}
f(t,k)\ :=\ |V\setminus\good{t}{k}|  
\end{equation}
the total number of non-$(t,k)$-good vertices. From the definition
of our path partitioning, we immediately obtain
\begin{equation}
  \label{eq:reachbound2}
  \expect{\length{\pi}}\ \le\ f(t,k)+t\cdot\expect{n_{t,k}(\pi)}\ .
\end{equation}

Adding up nonnegative multiples of~(\ref{eq:reachbound1})
and~(\ref{eq:reachbound2}) in such a way that
$\expect{n_{t,k}(\pi)}$ cancels out, one obtains the following bound:
$$
  \expect{\length{\pi}}\ \le\ 
  \frac{tn+\Tilde{g}(t,k)(f(t,k))}{\Tilde{g}(t,k)+t}\ \le\
  \frac{t}{\Tilde{g}(t,k)}n+f(t,k)
$$
This yields the following estimation.

\begin{lemma}
  \label{lem:reachgen}
  For $t,k\in\N$ with $t\ge 2$ and $k\ge 1$, we have
  $$
  \expect{\length{\pi}}\ \le\ 
  \frac{4td^t}{g(t,k)(g(t,k)-2t)}n+f(t,k)\ .
  $$
\end{lemma}


A general way to bound the function $f(t,k)$ is as follows.
\begin{lemma}
  \label{lem:bound-f}
  For $t,k\in\N$, we have
  $$
  f(t,k)\ \le\ \frac{d^t-1}{d-1}h_{<k}\ ,
  $$
  where $h_{<k}:=\sum_{j=0}^{k-1} h_j$ is the number of vertices with
  outdegree less than~$k$.
\end{lemma}

\begin{proof}
  If $v\in V\setminus\good{t}{k}$, then there is some
  $w\in\reach{t-1}{v}$ with $|\out{w}|<k$. On the other hand, each~$w$
  is contained in at most $\sum_{i=0}^{t-1}d^i=\tfrac{d^t-1}{d-1}$
  $(t-1)$-reaches (since the undirected graph is $d$-regular). 
  The claim follows.
\end{proof}

The following describes a way of bounding the function $g(t,k)$ by
studying the undirected graph of the polytope.

\begin{definition}[$(t,k)$-neighborhood, $\gamma(t,k)$]
  Let $t,k\in\N$. 
  \begin{enumerate}
  \item A subset $N\subset V$ is called a \emph{$(t,k)$-neighborhood}
    of~$v\in V$ if $N=\{v\}$ in case of $t=0$, or, if $t\ge 1$, there
    are~$k$ neighbors $w_1$, \dots, $w_k$ of~$v$ in the graph of~$P$
    together with $(t-1,k)$-neighborhoods $N_1$, \dots,
    $N_k$ of $w_1$, \dots, $w_k$, respectively, such that $N=\bigcup_{i=1}^k N_i$.
  \item We define $\gamma(t,k)$ as the minimum cardinality of
    $\setdef{w\in N}{\dist{v}{w}=t}$, taken over all $v\in V$ and all
    $(t,k)$-neighborhoods of~$v$. (Here, $\dist{v}{w}$ denotes the
    graph-theoretical distance between~$v$ and~$w$ in the undirected
    graph of~$P$.)
  \end{enumerate}
\end{definition}

If $v$ is $(t,k)$-good, then it follows 
right from the definitions that the boundary $\bdreach{t}{v}$ of its 
$t$-reach contains a $(t,k)$-neighborhood~$N$ of $v$.
In particular, all vertices $w\in N$ with $\dist{v}{w}=t$
are in $\bdreach{t}{v}$, and these are the ones of use to us.

\begin{lemma}
  \label{lem:bound-g}
  For $t,k\in\N$ with $t\ge 2$, we have
  $$
  g(t,k)\ \ge\ \gamma(t,k)\ . 
  $$
\end{lemma}

\subsection{Specialization to Cubes}
\label{subsec:cubes:spec}

In order to obtain from Lemma~\ref{lem:reachgen} an explicit bound for
the expected number of vertices visited by the \RE\ simplex-algorithm
on the $d$-cube, we will derive estimates on the functions $f(t,k)$
and $g(t,k)$ for $k=\lfloor\tfrac{d}{4}\rfloor$.

\begin{lemma}
  \label{lem:f-cubes}
  There is a constant $0<\alpha <1$
  such that
  $$
  f\big(t,\lfloor\tfrac{d}{4}\rfloor\big)\ \le\ 
  2^{\alpha d + \littleO{d}}
  $$
  holds for all $t\in\N$ 
  (where $f$ is the function defined in~(\ref{eq:def-f}) for the case
  of the $d$-cube, and with $k=\lfloor\tfrac{d}{4}\rfloor$).
\end{lemma}

\begin{proof}
  In the case of a $d$-cube and $k=\lfloor\tfrac{d}{4}\rfloor$, we have
  $$
  h_{<k}  
  \ =\ \sum_{i=0}^{\lfloor\tfrac{d}{4}\rfloor-1}\binom{d}{i}
  \ =\ 2^{h(\tfrac{1}{4})d+\littleO{d}}\ ,
  $$
  where $h(x)= x\log\frac{1}{x} + (1-x)\log\frac{1}{1-x}$ is the
  binary entropy function (see, e.g., \cite[Chap.~9, Ex.~42]{gkp-cm-94}).
  By Lemma~\ref{lem:bound-f} this implies the claimed bound (with the
  $\littleO{d}$ term depending on $t$).
\end{proof}

The final building block of our bound for the special case of cubes
is the following. Here, we denote by $\lofac{a}{b}$ (\emph{falling
factorial power}) the product $a(a-1)\cdots(a-b+1)$ (for $a,b\in\N$).

\begin{lemma}
  \label{lem:gamma-cubes}
  Let $t,k\in\N$ with $1\le t,k\le d$. If the polytope~$P$ considered in
  Section~\ref{subsec:cubes:gen} is a $d$-cube, then the following is
  true:
  \begin{enumerate}
  \item $\displaystyle
    \gamma(t,k)\ \ge\ \frac{k^t}{t!}-\sum_{i=1}^{t-1}\frac{k^i}{\lofac{t}{i}}\binom{d-1}{t-i-1}$.
  \item If~$t$ is a constant, then
    $\displaystyle
    \gamma(t,\lfloor\tfrac{d}{4}\rfloor)\ = \Omega(d^t)$.
  \end{enumerate}
\end{lemma}

\begin{proof}
  Part~(2) follows immediately from part~(1), since the sum becomes a
  polynomial in~$d$ of degree $t-1$ for $k=\tfrac{d}{4}$ (and constant~$t$).

  Let us prove~(1) for each fixed~$k$, by induction on~$t$, where the
  case $t=1$ holds due to $\gamma(1,k)=k$. Thus, let us consider the
  case $t\ge 2$.

  We may assume that the vertex~$v$ and its neighbors $w_1$, \dots,
  $w_k$, for which the minimum $\gamma(t,k)$ is attained, are
  $v=\allzero$ and $w_i=\evec{i}$ ($1\le i\le k$). For each~$i$, the
  $(t-1,k)$ neighborhood $N_i$ of~$\evec{i}$ has at least
  $\gamma(t-1,k)$ vertices~$w$ with $\dist{\evec{i}}{w}=t-1$, by
  definition. All of
  them have distance~$t-2$ or~$t$ from~$\allzero$. The former may be
  the case at most $\binom{d-1}{t-2}$ times (these vertices cannot
  have a one at position~$i$). Therefore, we have
  $$
  \big|\setdef{w\in N_i}{\dist{\allzero}{w}=t}\big|\ \ge\ \gamma(t-1,k)-\binom{d-1}{t-2}\ .
  $$

  On the other hand, every vertex~$w\in N_i$ with
  $\dist{\allzero}{w}=t$ needs to have a one at position~$i$
  (otherwise, $\dist{\evec{i}}{w}=t+1$). Hence, every vertex~$w$ with
  $\dist{\allzero}{w}=t$ can be contained in at most~$t$ of the
  neighborhoods $N_1$, \dots, $N_k$. Thus, we conclude (for $t\ge 2$)
  $$
    \gamma(t,k)\ \ge\ \frac{k\big(\gamma(t-1,k)-\binom{d-1}{t-2}\big)}{t}\ ,
  $$
  and thus, 
  \begin{equation}
    \label{eq:gamma:1}
    \gamma(t,k)\ \ge\ \frac{k}{t}\gamma(t-1,k)-\frac{k\binom{d-1}{t-2}}{t}\ .
  \end{equation}

  Using the induction hypothesis and (\ref{eq:gamma:1}) we derive
  $$
  \begin{array}{rcl}
    \gamma(t,k) \ 
    &\ge\ & \displaystyle
      \frac{k}{t}
      \left(
        \frac{k^{t-1}}{(t-1)!}-\sum_{i=1}^{t-2}\frac{k^i}{\lofac{(t-1)}{i}}\binom{d-1}{t-i-2}
      \right) 
      -\frac{k}{t}\binom{d-1}{t-2}\\
    &=\ & \displaystyle
      \frac{k^t}{t!}
      - \sum_{i=1}^{t-2}\frac{k^{i+1}}{\lofac{t}{i+1}}\binom{d-1}{t-i-2}
      - \frac{k^{0+1}}{\lofac{t}{0+1}}\binom{d-1}{t-0-2}\\
    &=\ & \displaystyle
      \frac{k^t}{t!}
      - \sum_{i=0}^{t-2}\frac{k^{i+1}}{\lofac{t}{i+1}}\binom{d-1}{t-i-2}\ ,
  \end{array}
  $$
which, after an index shift in the sum, yields the claim.
\end{proof}

Now we can prove our main result:
\begin{theorem}
  For every fixed $t\in\N$, there is a constant $C_t\in\R$ (depending
  on~$t$), such that the expected number of vertices visited by the \RE\
  simplex-algorithm on a $d$-dimensional cube, equipped with an
  abstract (in particular: a linear) objective function, is bounded by
  $$
  C_t\cdot\frac{2^d}{d^t}\ .
  $$
\end{theorem}

\begin{proof}
  Let $\pi$ be the (random) path (for some arbitrary start vertex)
  defined by the \RE\ simplex-algorithm on a $d$-cube equipped
  with an acyclic unique sink orientation. By
  Lemma~\ref{lem:reachgen}, we have, with $d':=\lfloor\tfrac{d}{4}\rfloor$,
  \begin{equation}
    \label{eq:thmcube:1}
  \expect{\length{\pi}}\ \le\ 
  \frac{4td^t}{g(t,d')(g(t,d')-2t)}2^d+f(t,d')\ .
  \end{equation}
  From Lemma~\ref{lem:f-cubes} we know that there is some constant 
  $0<\alpha<1$ with
  \begin{equation}
    \label{eq:thmcube:2}
    f(t,d')\ \le\ 2^{\alpha d + \littleO{d}}\ .
  \end{equation}
  Finally, by Lemmas~\ref{lem:bound-g} and~\ref{lem:gamma-cubes}~(2)
  there is some constant $\beta>0$ such that
  \begin{equation}
    \label{eq:thmcube:3}
    g(t,d')\ \ge\ \beta d^t\ .
  \end{equation}
  Putting~(\ref{eq:thmcube:1}), (\ref{eq:thmcube:2}),
  and~(\ref{eq:thmcube:3}) together, we obtain
  $$
    \expect{\length{\pi}}\ \le\ 
    \frac{4 t d^t}{\beta^2 d^{2t}-2t\beta d^t}2^d + 2^{\alpha d + \littleO{d}}\ ,
  $$
  which implies the claim.
\end{proof} 

\section{Conclusion}
\label{se:conclusion}

Probably
one can extend the methods we have used for analyzing \RE\ on
cubes to other classes of polytopes (e.g., general products of
simplices). However, it seems to us that it would be more interesting
to find a way of sharpening our bounds by enhancing our approach with some
new ideas.
As mentioned at the end of Section~\ref{sec:general}, the
analysis of our approach is sharp in the general setting. We suspect
that one cannot prove a subexponential bound for \RE\ on cubes with
our methods. Therefore, it would be most interesting to find a way of
combining our kind of analysis with some other ideas.

\bibliographystyle{plain}

\end{document}